\documentclass{article}
\usepackage{arxiv}

\usepackage[utf8]{inputenc} 
\usepackage[T1]{fontenc}    
\usepackage{hyperref}       
\usepackage{url}            
\usepackage{booktabs}       
\usepackage{amsfonts}       
\usepackage{nicefrac}       
\usepackage[english]{babel}
\usepackage{ae}
\usepackage{eucal}
\usepackage[dvips]{graphicx}
\usepackage{epsfig}
\usepackage{graphicx}
\usepackage{amsfonts,amsthm}
\usepackage{amssymb}
\usepackage{a4}
\usepackage{apu}
\usepackage{xcolor}
\usepackage{amsmath}
\usepackage{mathtools}
\usepackage{multicol}
\usepackage{subcaption}

\usepackage{dcolumn,amsthm}

\newtheoremstyle{theorem}{6pt}{6pt}{\rm}{}{\sffamily}{ }{ }{}
\theoremstyle{theorem}

\newtheoremstyle{lemma}{6pt}{6pt}{\rm}{}{\sffamily}{ }{ }{}
\theoremstyle{lemma}

\newtheoremstyle{example}{6pt}{6pt}{\rm}{}{\sffamily}{ }{ }{}
\theoremstyle{example}

\newtheoremstyle{corollary}{6pt}{6pt}{\rm}{}{\sffamily}{ }{ }{}
\theoremstyle{corollary}

\newtheoremstyle{definition}{6pt}{6pt}{\rm}{}{\sffamily}{ }{ }{}
\theoremstyle{definition}

\newtheoremstyle{remark}{6pt}{6pt}{\rm}{}{\sffamily}{ }{ }{}
\theoremstyle{remark}

\newtheoremstyle{approximation}{6pt}{6pt}{\rm}{}{\sffamily}{ }{ }{}
\theoremstyle{approximation}

\newtheoremstyle{scheme}{6pt}{6pt}{\rm}{}{\sffamily}{ }{ }{}
\theoremstyle{scheme}

\usepackage{csquotes}

\title{A Method for Finding Exact Solutions to the 2D and 3D Euler--Boussinesq Equations in Lagrangian Coordinates}

\author{
   Tomi Saleva \\
   Department of Physics and Mathematics\\
   University of Eastern Finland\\
   P.O. box 111, FI-80101 Joensuu, Finland  \\
   \text{tomisal@student.uef.fi} \\
      \And
   Jukka Tuomela \\
   Department of Physics and Mathematics\\
   University of Eastern Finland\\
   P.O. box 111, FI-80101 Joensuu, Finland  \\
   \text{jukka.tuomela@uef.fi} \\
}

\begin{document}

\maketitle

\begin{abstract}
We study the Boussinesq approximation for the incompressible Euler equations using Lagrangian description. The conditions for the Lagrangian fluid map are derived in this setting, and a general method is presented to find exact fluid flows in both the two-dimensional and the three-dimensional case. There is a vast amount of solutions obtainable with this method and we can only showcase a handful of interesting examples here, including a Gerstner type solution to the two-dimensional Euler--Boussinesq equations. In two earlier papers we used the same method to find exact Lagrangian solutions to the homogeneous Euler equations, and this paper serves as an example of how these same ideas can be extended to provide solutions also to related, more involved models.
\end{abstract}

\textbf{\emph Mathematics Subject Classification (2020)} 35A09; 35Q31; 76B70

\keywords{Euler equations, Boussinesq equations, Explicit solutions, Lagrangian formulation, Stratified fluids}

\thanks{The first author was supported by Finnish Cultural Foundation.}

\section{Introduction}
We continue our investigations of  finding explicit solutions to the incompressible Euler equations in the Lagrangian framework. In  \cite{toju,toju2} we considered the homogeneous case while here we treat the Boussinesq approximation of the heterogeneous case. The idea behind the Boussinesq approximation is that the density of the fluid does not fluctuate very much around the mean value, so that one supposes that the density is constant, except in the term that represents the gravity. There are many physical situations where this simplification is justified. The general introduction to this type of models can be found in \cite{majda}, while in regard to the Lagrangian description in general we refer to \cite{bennett}.

One of the first explicit solutions to the homogeneous Euler equations in Lagrangian formulation is due to Gerstner \cite{G}, and it turns out that solutions of the Gerstner type show up in many related models.
For example the Gerstner solution plays a prominent role also in the heterogeneous Euler equations, even without using the Boussinesq approximation.  In \cite{dubreil} it was shown that the Gerstner wave can also be a barotropic flow, i.e. a flow where pressure and density are functions of each other, see also a modern exposition by Stuhlmeier \cite{Stuhlmeier2011}. Gerstner type solutions are also relevant in many physical models that are somehow modifications of the Euler equations. For example Constantin \cite{con4} found a barotropic Gerstner type exact solution to the equatorial water wave equations with the beta-plane approximation of the Coriolis effect. For more papers on the applicability of the Gerstner waves see for example \cite{henry,conmon,con1,con2} and especially the survey articles \cite{AP,Stuhlmeier2015}. In \cite{Weber} there is also a Gerstner type solution to the first and second terms of the asymptotic expansion to the heterogeneous Euler equations. One remarkable property of Gerstner type solutions that makes them popular, is that they can in many cases be interpreted as free surface solutions, or in other words one can model the interface of two different fluids. Apparently there are no other known explicit solutions with this property.

Many explicit solutions to stratified fluid flows are shear flows, which are flows with two-dimensional horizontal motion that only depends on height. Most of the studies use the Eulerian framework. As for the Lagrangian description, Yakubovich and Shrira \cite{YS} found solutions with columnar motion using the Boussinesq approximation, but it seems that elsewhere the Eulerian framework is used.

As stated above there are numerous models for geophysical fluid flows, depending on the particular context. We have chosen to analyze the Euler--Boussinesq model in the present article. For a thorough discussion of this model and its applicability we refer to \cite{majda}.  As in \cite{toju,toju2}, we use a separation of variables method to compute our solutions, that is, we search for fluid particle maps that can be expressed as the product of a time-only dependent matrix and a spatial-only dependent vector. This approach leads for example to a two-dimensional Gerstner type flow \eqref{24ell}, as well as a plethora of other solutions both in the 2D and the 3D setting. Solution \eqref{24ell} does not satisfy the free boundary condition like most Gerstner type flows, but it still gives valid internal flows, which are an important application of the Euler--Boussinesq equations. Note in particular that our approach is not restricted to the Euler equations and Euler--Boussinesq equations, and we are confident that it could be used profitably to analyze other related models.  For example Abrashkin \cite{abrashkin2019} derived the governing Lagrangian equations for equatorial beta-plane flows, allowing anyone to readily attempt the separation of variables strategy that we use below. Also the Gerstner type solutions in all those different models discussed above could have been found by our method.

This paper is organized as follows. In Section 2 we give basic definitions and specify the model which will be analyzed. In Section 3 we discuss columnar flows and generalize one family of solutions which was given in \cite{YS}. In Section 4 we find exact solutions to the two-dimensional Euler--Boussinesq equations while utilizing the similarities to the homogeneous situation \cite{toju}. In Section 5 we consider the three-dimensional case and use the framework outlined in \cite{toju2} to find solutions to the Euler--Boussinesq equations. 

In some cases we are able to find explicit periodic and nonperiodic solutions. In other cases one can show that a solution exists for all times for convenient parameter values. Finally there are cases where all one can say that a local solution is well-defined. We allow for both stably and unstably stratified fluids as both types of situations can often be covered by a single formula, though when we give specific examples we concentrate on the stably stratified case. Some solutions that we find are neither stably nor unstably stratified. We also note that it has not been possible to explore all different cases in the present article.

\section{Euler--Boussinesq equations}

\subsection{Notation}
Let $A \,:\,\mathbb{R}\to\mathbb{R}^{n\times m}$ and $v\,:\, \mathbb{R}^n\to \mathbb{R}^m$. 
We denote the columns of $A$ by $A_i$ and its entries by $a_{ij}$. Now the minors of $A$ are denoted by
\begin{equation}
    \begin{aligned}
        p_{ij}=&\det(A_i,A_j)\ ,\,
        \mathrm{if}\ n=2\ ,\\
         p_{ijk}=&\det(A_i,A_j,A_k)\ ,\,
        \mathrm{if}\ n=3\ .
    \end{aligned}
    \label{min-A}
\end{equation}
Similarly the the minors of $dv$ are
\begin{equation}
    \begin{aligned}
        g_{ij}=&\det(\nabla v^i,\nabla v^j)\ ,\,
        \mathrm{if}\ n=2\ ,\\
         g_{ijk}=&\det(\nabla v^i,\nabla v^j,\nabla v^k)\ ,\,
        \mathrm{if}\ n=3\ .
    \end{aligned}
    \label{min-dv}
\end{equation}
It will also be useful to define
\begin{equation}
    \begin{aligned}
    Q_{ij}=&\langle A_i',A_j\rangle - \langle A_j',A_i\rangle \ ,\\
        G_{ij}=&\nabla v^i \times \nabla v^j\ ,
    \end{aligned}
    \label{QjaG}
\end{equation}
where $A_i'$ is the (time) derivative of $A_i$ and $G_{ij}$ is useful only when $n=3$. For derivatives of $v$ we use multiindices so that for example 
\[
   v^1_{201}=\frac{\partial^3v^1}{\partial z_1^2\partial z_3}\ .
\]
At times we will also say that functions $v^1$ and $v^2$ are an anti-Cauchy--Riemann pair, or an anti-CR pair, if
\[
\begin{cases}
    v^1_{10}+v^2_{01}=0\\
    v^1_{01}-v^2_{10}=0
\end{cases}\ .
\]

\subsection{Model}
The $n$-dimensional heterogeneous incompressible Euler equations are given by the system
\begin{equation}
\begin{aligned}
    \nabla\cdot u&=0 \\
    \tilde \rho (u_t+u\nabla u+ g e_n)+\nabla p&=0  \\
    \tilde \rho_t+\langle u,\nabla \tilde \rho\rangle &= 0\ .
\end{aligned}
\label{euler}
\end{equation}
Here $u$ is the velocity field, $\tilde \rho$ is the density, $g$ is the acceleration due to gravity, and $e_n$ is the vertical unit vector parallel to gravity.
Let us write these equations in the Lagrangian framework. Let $D\subset\mathbb{R}^n$ be a domain and let us consider a family of diffeomorphisms $\varphi^t\,:\,D\to \Omega_t=\varphi^t(D)$. The coordinates in $D$ are denoted by  $z$ and in $\Omega_t$ by $x$. We can also define 
\[
  \varphi\,:\,D\times \mathbb{R}\to \mathbb{R}^n\quad,\quad
  \varphi(z,t)=\varphi^t(z)\, .
\]
Now given such $\varphi$ we can define the associated vector field $u$ by the formula
\begin{equation}
\frac{\partial}{\partial t} \fii(z,t)=u(\fii(z,t),t)\,.
\label{siirto}    
\end{equation}
Then $(u,\tilde \rho,p)$ solves \eqref{euler} if $\det(d\varphi)\ne 0$ and 
\begin{subequations}
\begin{align}
    \frac{d}{dt}\det(d\fii)&=0 
\label{lagrange1}\\
    \hat \rho \big(d\fii^T(\fii''+ g e_n)\big)+\nabla \hat p &= 0
\label{lagrange2}\\
    \frac{d}{dt} \hat \rho &= 0\ .
\label{lagrange3}
\end{align}
\end{subequations}
where $\hat p=p\circ \varphi$ and $\hat \rho=\tilde \rho\circ \varphi$. Hence $\hat \rho$ is  a function of spatial variables only, which is one great advantage of using the Lagrangian description.\footnote{Strictly speaking in the Lagrangian description we should have $\det(d\varphi)=1$. However, given $\varphi$ as above we can define $\Phi^t=\varphi^t\circ (\varphi^0)^{-1}$.} Typically one cannot explicitly recover $u$ from $\fii$, since it requires computing the inverse of $\fii^t$. One exception to this is obtained when $\fii=A(t)z$ for some square matrix $A$, yielding $u=A'A^{-1}x$.

The standard way to apply the Boussinesq approximation is to assume that $\hat \rho$ is constant in every term except when it is multiplied by $g$. Thus the Boussinesq approximation takes into account how density variations affect buoyancy. This would mean that Newton's second law \eqref{lagrange2} would be replaced by the equation
\begin{equation}
d\fii^T(\overline \rho \fii''+ \hat \rho g e_n)+\nabla \hat p=0\ ,
\label{newton2}
\end{equation}
where $\overline \rho$ is the average density. However, as was shown in \cite{YS}, we can make the model slightly more accurate without making the equations any more difficult to study. Supposing that $n=3$ and taking the curl of \eqref{lagrange2}, we obtain
\begin{equation}
\nabla \hat \rho \times (d\fii^T\fii'') + \hat \rho \Big(\sum_{i=1}^3 \nabla \fii_i'' \times \nabla \fii_i \Big) + g\nabla \hat \rho \times \nabla \fii_3 = 0\,.
\label{newton}
\end{equation}
Now assuming only $\nabla \hat \rho \times (d\fii^T\fii'')=0$ in \eqref{newton}, dividing  by $\hat \rho$ and defining $\rho = g\ln{\hat \rho} $
we obtain
\begin{equation}
    \sum_{i=1}^3 \nabla \fii_i'' \times \nabla \fii_i + \nabla \rho \times \nabla \fii_3 = 0\ .
\label{newton3}
\end{equation}
If we had used \eqref{newton2}, we would have obtained the same equation but with $\rho$ defined as $\rho=g\hat \rho/\overline \rho$ instead of $\rho = g\ln{\hat \rho} $. Hence it should be kept in mind that $\rho$ is not really density, but the density $\hat\rho$ can be recovered from $\rho$ using either $\rho=g\hat \rho/\overline \rho$ or $\rho = g\ln{\hat \rho} $.

Integrating the left hand side of \eqref{newton3} with respect to $t$ we obtain the equivalents of what are the Cauchy invariants for the Euler equations:
\begin{equation}
    h=(h^1,h^2,h^3) = \sum_{i=1}^3 \nabla \fii_i' \times \nabla \fii_i + \nabla \rho \times \nabla \int \fii_3\textrm{ }dt\,.
\label{cauchy}
\end{equation}
This integrated form is often useful since it removes the second-order dependence of $\fii_1$ and $\fii_2$. But calling the components of $h$ the "Cauchy invariants of the Euler--Boussinesq equations" is perhaps a bit of a stretch since the time integral of $\fii_3$ produces an arbitrary function of $z$ and so there is no canonical way to choose $h$.

In the two-dimensional case $h$ is just a scalar and it is convenient to write it in the following form:
\begin{equation}
    h = \sum_{i=1}^2 \det\big(\nabla \fii_i' , \nabla \fii_i\big) + \det\Big(\nabla \rho , \nabla \int \fii_2\textrm{ }dt\Big)\, .
\label{cauchy2d}
\end{equation}
We have now established the conditions for a solution to the 2D and 3D Euler--Boussinesq equations:
\begin{theorem}
The pair $(\fii, \rho)$ provides a solution to the Euler--Boussinesq equations \eqref{euler} via \eqref{siirto} if and only if $\det(d\fii)\neq 0$ everywhere, and $\det(d\fii)$, $\rho$ and $h$ are independent of time, where $h$ is given by \eqref{cauchy2d} in the 2D case and by \eqref{cauchy} in the 3D case.
\label{yleinen}
\end{theorem}

\subsection{Separation of variables}

As in our previous papers \cite{toju,toju2}, we try to find solutions of the form $\fii(z,t)=A(t)v(z)$, where $A\,:\,\mathbb{R}\to \mathbb{R}^{n\times m}$ and $v\,:\,D\subset \mathbb{R}^n\to \mathbb{R}^m$. To this end we need convenient formulas for $\det(d\varphi)$ and $h$. Using \eqref{min-A}, \eqref{min-dv} and  the Cauchy--Binet formula we obtain
\begin{equation}
\begin{aligned}
    \det(d\fii)=&\sum_{1\leq i<j\leq m} p_{ij}g_{ij}
    \ ,\,
        \mathrm{if}\ n=2\ ,\\
        \det(d\fii) = &\sum_{1\leq i<j<k\leq m}p_{ijk}g_{ijk}
        \ ,\,
        \mathrm{if}\ n=3\ .
    \end{aligned}
    \label{det2dja3d}
\end{equation}
Then using \eqref{QjaG} we compute that
\begin{equation}
    \begin{aligned}
   h=&\sum_{1\leq i<j\leq m} Q_{ij}g_{ij} + \sum_{i=1}^m y_i\det(\nabla \rho,\nabla v^i)
    \ ,\,
        \mathrm{if}\ n=2\ ,\\
   h =& \sum_{1\leq i<j\leq m}Q_{ij}G_{ij} + 
    \sum_{i=1}^m y_i\nabla \rho \times \nabla v^i\
        \ ,\,
        \mathrm{if}\ n=3\ ,
    \end{aligned}
    \label{h2dja3d}
\end{equation}
where $y_i$ is a function such that $y_i'=a_{ni}$.

Some properties of $\fii$ and $\rho$ are gathered in the following Lemma.

\begin{lemma}
Let $(\fii,\rho)$ be a solution to the Euler--Boussinesq equations.
\begin{enumerate}
\item Let  $\psi\,:\,\hat D\to D$ be an arbitrary diffeomorphism and let $\tilde\fii^t=\fii^t\circ\psi$, $\rho_0 = \rho \circ \psi$. Then $(\tilde\fii,\rho_0)$  is a solution to the Euler--Boussinesq equations if and only if $(\fii, \rho)$ is.

\item Let $\fii=Av$. If $H$ is a regular $m\times m$ matrix with constant entries and $\tilde v=Hv$, $\tilde A = AH^{-1}$, then $(\tilde A \tilde v, \rho)$ is a solution.

\item In the 3D case, if $R$ is a constant rotation matrix such that $(R\fii)_3=\fii_3$, then $(R\fii,\rho)$ is also a solution. In the 2D case there is no nontrivial rotation $R$ for which $R\fii$ is always a solution.
\end{enumerate}
\label{properties}
\end{lemma}

The change of coordinates in the first part of this Lemma typically allows us to assume without loss of generality that $(v^1,v^2)=(z_1,z_2)$ in the 2D case and $(v^1,v^2,v^3)=(z_1,z_2,z_3)$ in the 3D case. In the present paper we always use this simplification, but in practice the general form allows for more flexibility since the inverse map needed for this transformation cannot always be explicitly computed. The second part of the Lemma can be used to bring the problems to simpler form without loss of generality. From this part it also follows that if $A_i$ or $\nabla v^i$ are linearly dependent over $\mathbb{R}$, the solution reduces to a case of lower $m$. On the other hand, $\nabla\rho$ being an $\mathbb{R}$-linear combination of $\nabla v^i$ does not imply that the solution reduces in this way.

The formulas \eqref{det2dja3d} and \eqref{h2dja3d} allow us to deduce what kind of constraints we should set for the spatial functions $v^i$ and $\rho$. We want the time derivatives of $\det(d\fii)$ and $h$ to vanish for all $t$, so, for example in the two-dimensional case, fixing any $t$  produces constraints of the form
\begin{align*}
\sum_{1\leq i<j\leq m} \beta_{ij}g_{ij}&=0\ , \\
\sum_{1\leq i<j\leq m} \gamma_{ij}g_{ij} + \sum_{i=1}^m \gamma_i\det(\nabla \rho,\nabla v^i)&=0\ ,
\end{align*}
where $\beta_{ij}$, $\gamma_{ij}$ and $\gamma_i$ are constants. Thus the spatial functions need to satisfy a number of constraints like these. By substituting these constraint equations back to the formulas of $\det(d\fii)$ and $h$ we also immediately obtain the conditions that the time functions $a_{ij}$ are required to satisfy. In \cite{toju} and \cite{toju2} we derived the spatial constraints from the formulas of $\det(d\fii)$ and $h$ for the homogeneous Euler equations, and especially in \cite{toju} we went into greater detail in the 2D case to show what were all possible constraint sets that are essentially different for the cases that we studied. In the case of the Euler--Boussinesq equations this analysis, which is done using the second part of Lemma \ref{properties}, yields very similar results, so in the present paper we mainly consider cases that are similar to those in \cite{toju} and \cite{toju2}. Here the presence of $\rho$ in the formula of $h$ makes the analysis slightly different, and while we consider several different possibilities for $\rho$ as well, we omit large amount of situations that where the interplay between $\rho$ and $v$ is more intricate.

The equations for $h$ are second-order differential equations for $a_{nj}$. This makes the problem very hard in general, but we can still find some exact solutions, though we often need to restrict to special cases where there are not many terms in $\fii_n$. Even so, the number of different solutions we can find is so vast that we have to restrict to the cases that seem interesting as well as  sufficiently simple. Many times, when we have an underdetermined system and are able to choose some functions arbitrarily, we choose the ones that have these second-order derivatives. Thus the system becomes a first-order system in terms of the unknown functions and is more easily solvable. However, the resulting solution formulas are sometimes quite ugly, and it is possible to analyze the systems with other approaches as well.

\section{Columnar flows}
Before studying any specific cases, we would like to show a general method of obtaining solutions to the Euler--Boussinesq equations from a specific type of flows that satisfy the homogeneous Euler equations. We consider the 3D case first, the 2D case is similar.

In many cases that we considered in \cite{toju2}, we could find solutions to the 3D Euler equations that were of the form
\begin{equation}
\fii=\big(\fii_1(z_1,z_2,t),\fii_2(z_1,z_2,t),a(t)z_3\big)\ ,
\label{columnar}
\end{equation}
where $a$ is always nonzero. Flows like this feature columnar motion, where the vertical columns can stretch and contract but otherwise remain intact as the flow evolves. These are also solutions to the Euler--Boussinesq equations when density is an arbitrary function of $z_3$, as $\nabla \rho \times \nabla \fii_3$ vanishes and condition \eqref{cauchy} reduces to the usual Cauchy invariants of the Euler equations. 
But in this case we can look for more solutions to the Euler--Boussinesq equations by choosing $\rho=c_0z_3$ and
\[
\fii=\big(\fii_1(z_1,z_2,t),\fii_2(z_1,z_2,t),f(z_1,z_2,t)+a(t)z_3\big)\ .
\]
In \cite{YS} a particular solution of this form was presented in the stably stratified special case $a(t)=1$, $c_0<0$. In this case  the general solution of $\fii_3$ is
\begin{equation}
\fii_3 = z_3+f^1(z_1,z_2)\cos(Nt)+f^2(z_1,z_2)\sin(Nt)\ ,
\label{YS}
\end{equation}
where the constant $N$, given by $N^2=-c_0$, is the Brunt--V\"ais\"al\"a frequency, $f^i$ are arbitrary, and $\fii_1$ and $\fii_2$ have to satisfy the usual 2D Euler equations. 

Actually one can describe the solutions for arbitrary $a$. First we note that $f$ has no effect on $\det(d\fii)$ and $h^3$, whereas for the time derivatives of $h^1$ and $h^2$ we have
\begin{align*}
    (h^1)'&=f_{010}''a-a''f_{010}-c_0 f_{010} \\
    -(h^2)'&=f_{100}''a-a''f_{100}-c_0f_{100}\ .
\end{align*}
Both $(h^1)'$ and $(h^2)'$ vanish if and only if
\begin{equation}
af''-(a''+c_0)f = 0\ .
\label{moi}
\end{equation}
This is a linear ODE where $z$ appears only as a parameter so the general solution is of the form 
\[
f=a_1(t)f^1(z_1,z_2)+a_2(t)f^2(z_1,z_2)\ .
\]
Substituting this back to \eqref{moi} we see that $a_i$ are two linearly independent solutions of 
\begin{equation}
y''+qy=0\ \mathrm{, where}\  q=-(a''+c_0)/a\,.
\label{q-yht}
\end{equation}
Note that $q$ and hence $a_i$ are well-defined since we must anyway choose $a$ such that $a\ne 0$ for all $t$.

So the solution set is parametrized by $a$ via $q$, and now by classical theorems if we choose $a$ such that  $q>0$ then $a_j$ are  oscillating and bounded solutions.
Note that there is no condition whatsoever for $f^1$ and $f^2$, other than that they do not depend on $z_3$. 

Of course we can also choose $c_0=0$ to obtain new solutions to the usual Euler equations, although in that case $a$, $a_1$, and $a_2$ are linearly dependent and we may assume $a_2=0$.

In the 2D situation choosing $\rho=c_0z_2$   the analogous  $\fii$ is given by
\begin{equation}  \fii=\big(\fii_1(z_1,t),f(z_1,t)+a(t)z_2\big)\ .
    \label{columnar2d}
\end{equation}
The incompressibility condition \eqref{lagrange1} already says that $\det(d\fii)=(\fii_1)_{10}a(t)$ is independent of time, so we may assume by a coordinate transformation that $\fii_1=z_1/a(t)$.
Then the formula of $h$ gives us the same condition as in the three-dimensional case so the solution can be written as
\begin{equation}
\fii=\big(z_1/a,a_1f^1(z_1)+a_2f^2(z_1)+a z_2\big)
  =Av=
  \begin{pmatrix}
      1/a&0&0&0\\
      0&a&a_1&a_2
  \end{pmatrix}
  \begin{pmatrix}
      z_1\\z_2\\f^1\\f^2
  \end{pmatrix}\ ,
  \label{2d-col}
\end{equation}
where $a_i$ are the solutions of \eqref{q-yht} and $f^i$ are arbitrary. This is the only 2D solution of this form, so below we will only use the generalization method of this Section in 3-dimensional situations.

\section{2-dimensional case}
\subsection{$m=2$}
Let us start by studying the two-dimensional case, which is the easier case. In case $m=2$ we always have $v=(z_1,z_2)$ by part 1 of Lemma \ref{properties}, and by \eqref{det2dja3d} and \eqref{h2dja3d} $A$ satisfies
\begin{equation}
\begin{aligned}
    \det(d\fii)&=a_{11}a_{22}-a_{12}a_{21} = 1 \\
    h&=a_{11}'a_{12}-a_{12}'a_{11}+a_{21}'a_{22}-a_{22}'a_{21}-y_1\rho_{01}+y_2\rho_{10}\ .
\end{aligned}
\label{2x2}
\end{equation}
Here $\det(d\fii)=1$ can be assumed by scaling $A$. By a linear transformation described in part 2 of Lemma \ref{properties} we may also assume that either $a_{21}=0$, or $a_{21}$ and $a_{22}$ are linearly independent.

If $a_{21}=0$, the incompressibility condition requires $a_{22}\neq 0$, and in this case we must have $\rho=c_0z_1+f(z_2)$ where $f$ is arbitrary. Then the equations are
\begin{align*}
    a_{11}a_{22}&=1 \\
    a_{11}'a_{12}-a_{12}'a_{11}+c_0y_2&=c\ ,
\end{align*}
where $c$ is a constant. If $a_{22}$ is chosen arbitrarily, then $a_{11}=1/a_{22}$ and $a_{12}=\frac{1}{a_{22}}\int (c_0y_2-c)a_{22}^2$ $dt$. We can also derive the Eulerian description for this solution by computing
\begin{align*}
\rho \circ \fii^{-1} &= c_0\Big(a_{22}x_1-(x_2/a_{22})\int(c_0y_2-c)a_{22}^2\textrm{ }dt\Big)+f(x_2/a_{22}) \ ,\\
u&=A'A^{-1}x=\begin{pmatrix}-a_{22}'/a_{22} & c_0y_2-c \\
0 & a_{22}'/a_{22}
\end{pmatrix}x\ .
\end{align*}

If $a_{21}$ and $a_{22}$ are both nonzero and linearly independent, then $\rho$ must be linear. In this case we may assume by a linear transformation that $\rho=c_0z_2$. If we write $A$ as
\[
A=\begin{pmatrix} 
\cos(\theta) & -\sin(\theta) \\
\sin(\theta) & \cos(\theta)
\end{pmatrix} \begin{pmatrix}
b & \ell b \\
0 & 1/b
\end{pmatrix}\ ,
\]
then the incompressibility condition is satisfied and the equation for $h$ gives
\[
b^2\ell'=2\theta' - c - c_0y_1\ ,
\]
where now $y_1'=a_{21}=b\sin(\theta)$. Thus we can choose $b$ and $\theta$ arbitrarily and solve for $\ell$ to obtain all the solutions. In Eulerian coordinates the solution is
\begin{align*}
\rho \circ \fii^{-1}&=
c_0b\big(\cos(\theta)x_2 - \sin(\theta)x_1\big) \ ,\\
u&=\begin{pmatrix}
    \cos(2\theta) & -\sin(2\theta) \\
    \sin(2\theta) & \cos(2\theta)
\end{pmatrix}\begin{pmatrix}
     b'/b & \theta' \\
     \theta' & -b'/b
\end{pmatrix}x+(c_0y_1+c)\begin{pmatrix}
   \sin(\theta) \cos(\theta) & -\cos^2(\theta) \\
    \sin^2(\theta) & -\sin(\theta)\cos(\theta)
    \end{pmatrix}
x\ .
\end{align*}

\subsection{$m=4$}
We turn immediately to case $m=4$, skipping the analysis of case $m=3$ completely. There are solutions in case $m=3$ for several choices of $\rho$ but there is no room for us to consider them here.

We meticulously showed in \cite{toju} (Theorems 4.1 and 4.2) for the homogeneous Euler equations that when $A$ is a $2\times 4$ matrix, then without loss of generality the spatial component $v$ may be assumed to satisfy one of these four cases:
\begin{enumerate}
    \item $g_{13}+g_{24}=g_{14}-g_{23}=0$,
    \item $g_{13}=g_{24}=0$,
    \item $g_{13}+g_{24}=g_{14}=0$, or
    \item $g_{14}=g_{24}=0$.
\end{enumerate}
The same is true for the Euler--Boussinesq equations with an identical proof since the incompressibility condition \eqref{det2dja3d}, which is the same for both systems, is all that is required to prove this claim.

We will consider all four cases in this Section. The equations for the entries of $A$ are not underdetermined in the first three cases so the second-order ODEs will make them difficult to study. We still find one special solution to all these cases. The fourth case, however, is underdetermined and we are able to present a general solution for one choice of $\rho$.

\subsubsection{Case 1}
Let $v$ satisfy $g_{13}+g_{24}=g_{14}-g_{23}=0$, in which case $v=(z_1,z_2,f^1,f^2)$, where $f^1$ and $f^2$ are an anti-CR pair. It is shown in \cite{toju} how we can use linear transformations to bring $A$ to a simpler form without loss of generality. In this case, \cite[Lemma 5.3]{toju} implies that we may assume that
\[
A=\begin{pmatrix} 
\cos(\theta) & -\sin(\theta) \\
\sin(\theta) & \cos(\theta)
\end{pmatrix} \begin{pmatrix}
b_{11} & b_{12} & \cos(\mu)b_{11} + \sin(\mu)b_{12} & \cos(\mu)b_{12} - \sin(\mu)b_{11} \\
0 & 1/b_{11} & \sin(\mu)/b_{11} & \cos(\mu)/b_{11}
\end{pmatrix}\ .
\]
In this form $A$ satisfies the incompressibility condition \eqref{det2dja3d} with
\[
  \det(d\fii) = 1- \lvert \nabla f^1\rvert ^2 \ ,
\]
which has to be nonzero in $D$.

Here $\rho$ could be chosen to be for example an arbitrary linear combination of $v^i$ but the equations would become too difficult for us to find any exact solutions. Instead we choose the more simple $\rho=c_0z_2$, in which case \eqref{h2dja3d} yields
\begin{align*}
h=&(Q_{12}-c_0y_1)g_{12}+Q_{34}g_{34}+(Q_{13}-Q_{24}-c_0y_4)g_{13}+(Q_{14}+Q_{23}+c_0y_3)g_{14} \\
=&c_{12}g_{12}+c_{34}g_{34}+c_{13}g_{13}+c_{14}g_{14}\ .
\end{align*}
We have chosen the subscripts of the constants $c_{ij}$ according to the corresponding $g_{ij}$. In later cases we often use this notation without further notice.

\begin{figure}
\centering
\includegraphics[width=\linewidth]{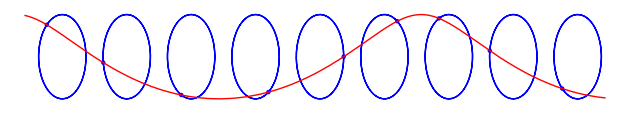}
\caption{Example of solution \eqref{24ell}: some ellipse-shaped trajectories and a curve of constant density at a fixed time. In this example $f^1=e^{z_2}\cos(z_1)$ and $f^2=e^{z_2}\sin(z_1)$.}
\label{ellipsit}
\end{figure}

The conditions we collect from the formula of $h$ are explicitly written as
\begin{align*}
    Q_{12}-c_0y_1 = 2\theta' + b_{11}'b_{12}-b_{12}'b_{11}-c_0y_1 &= c_{12} \\
    Q_{34} = 2\theta' + b_{11}'b_{12}-b_{12}'b_{11}+\mu'(b_{11}^2+b_{12}^2+1/b_{11}^2) &= c_{34} \\
    Q_{13}-Q_{24}-c_0y_4 = \mu'\sin(\mu)(b_{11}^2-b_{12}^2-1/b_{11}^2)-2\mu'\cos(\mu)b_{11}b_{12}-c_0y_4 &= c_{13} \\
    Q_{14}+Q_{23}+c_0y_3 = \mu'\cos(\mu)(b_{11}^2-b_{12}^2-1/b_{11}^2)+2\mu'\sin(\mu)b_{11}b_{12}+c_0y_3 &= c_{14}\ .
\end{align*}
We can bring this system to the form
\begin{align*}
b_{11}' =& -\big(c_0b_{11}b_{12}\cos(\theta) + 2b_{12}s^2 + c_0\sin(\theta)\big)/\big(4s\big) \\
b_{12}' =& \big(  c_0b_{11}^2b_{12}\sin(\theta)-c_0b_{11}^3b_{12}^2\cos(\theta) - 2c_0b_{11}\cos(\theta) + 2(b_{11}^4 - 1)s^2\big)/\big(4b_{11}^3s\big) \\
s' =& c_0\big(b_{11}b_{12}\cos(\theta) - \sin(\theta)\big)/(2b_{11}) \\
\theta_0'=&\Big(2c_0^2b_{11}^2b_{12}\sin^2(\theta) + c_0^2b_{11}^2b_{12} - 8b_{12}s^4 + \big(4c_0b_{11}^3b_{12}s\theta_0 - 10c_0b_{11}b_{12}s^2\big)\cos(\theta) \\
&+ \big(4c_0b_{11}^2s\theta_0 + 2(3c_0b_{11}^4 - c_0)s^2 - (3c_0^2b_{11}^3b_{12}^2 + 5c_0^2b_{11})\cos(\theta)\big)\sin(\theta)\Big)/\big(16b_{11}^3s^2\big)\\
\theta'=&\theta_0\ ,
\end{align*}
where $s=\mu'$. The solution blows up if $b_{11}$ or $s$ reaches zero, and it is difficult to analyze with which initial data this happens. However, the equilibrium point
\[
 \big(b_{11},b_{12},s,\theta_0,\theta\big)=
 \big(c_1,0,\mu_0,0,0   \big)\quad\mathrm{, where}\quad
 \mu_0^2=\frac{c_0c_1}{c_1^4-1}\ ,
\]
leads to a simple solution.
We may assume that $c_1>0$. Thus, if $c_0>0$ then $c_1>1$, and if $c_0<0$ then $c_1<1$. This gives 
\begin{equation}
A=\begin{pmatrix}
    c_1 & 0 \\
    0 & 1/c_1
\end{pmatrix}\begin{pmatrix}
    1 & 0 & \cos(\mu_0t) & -\sin(\mu_0t) \\
    0 & 1 & \sin(\mu_0t) & \cos(\mu_0t)
\end{pmatrix}\ ,
\label{24ell}
\end{equation}
This looks like a usual Gerstner type solution from the case of Euler equations \cite[Theorem 5.1]{toju}, except that $\fii_1$ and $\fii_2$ have been scaled. In case $c_0<0$ and $c_1<1$, which describes the stably stratified situation, the particle trajectories are ellipses stretched in the vertical direction. The curves of constant density are formed by the particles whose centers of trajectory are on the same horizontal line, see Figure \ref{ellipsit}.

\begin{figure}
\centering
\begin{subfigure}{0.15\linewidth}
\includegraphics[width=\linewidth]{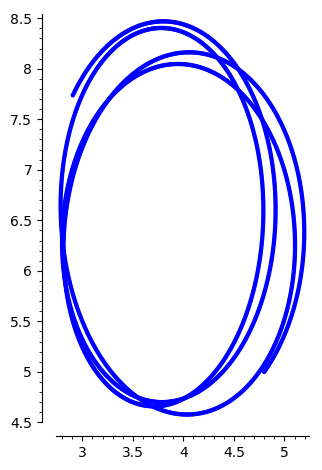}
\caption{$\delta=-0.02$}
\end{subfigure}
\begin{subfigure}{0.15\linewidth}
\includegraphics[width=\linewidth]{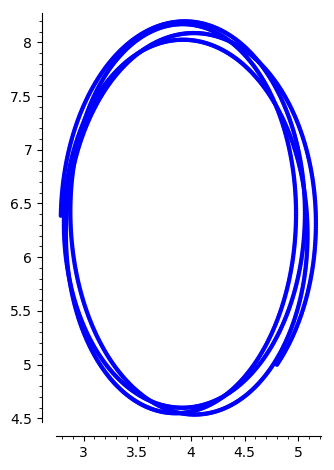}
\caption{$\delta=-0.01$}
\end{subfigure}
\begin{subfigure}{0.15\linewidth}
\includegraphics[width=\linewidth]{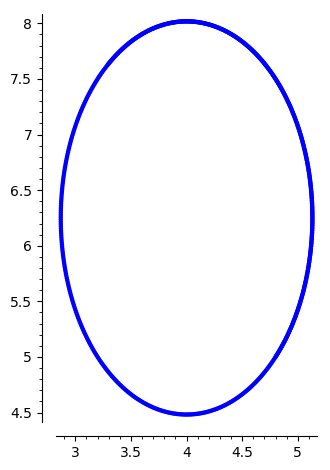}
\caption{$\delta=0$}
\end{subfigure}
\begin{subfigure}{0.15\linewidth}
\includegraphics[width=\linewidth]{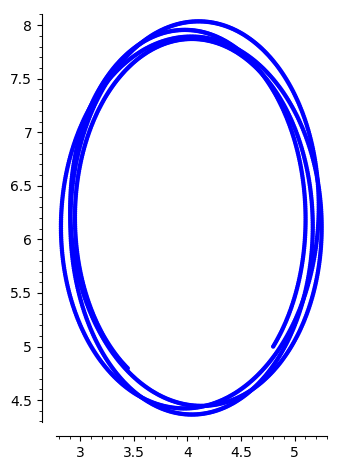}
\caption{$\delta=0.01$}
\end{subfigure}
\begin{subfigure}{0.15\linewidth}
\includegraphics[width=\linewidth]{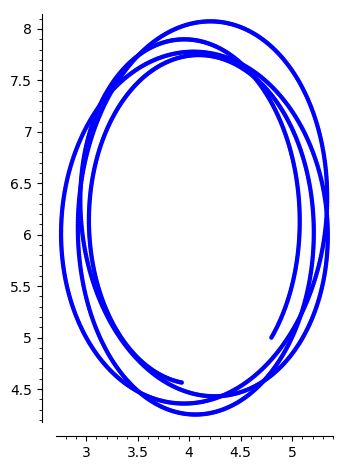}
\caption{$\delta=0.02$}
\end{subfigure}
\caption{Numerically solved trajectories for the general case where $c_0=-1$ and initially we have set $(b_{11},b_{12},s,\theta_0,\theta)=(0.8,0,\mu_0+\delta,0,0)$, where $\delta$ varies between $-0.02\dots 0.02$.}
\label{24radat}
\vspace*{\floatsep}
\begin{subfigure}[b]{0.15\linewidth}
\includegraphics[width=\linewidth]{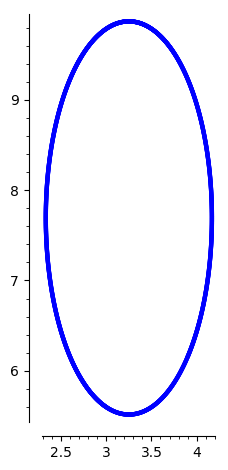}
\caption{$\delta=0$}
\end{subfigure}
\begin{subfigure}[b]{0.236\linewidth}
\includegraphics[width=\linewidth]{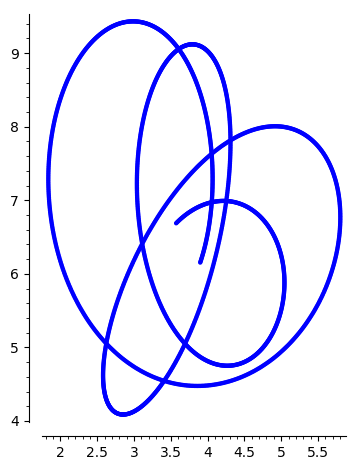}
\caption{$\delta=0.01$}
\end{subfigure}
\begin{subfigure}[b]{0.277\linewidth}
\includegraphics[width=\linewidth]{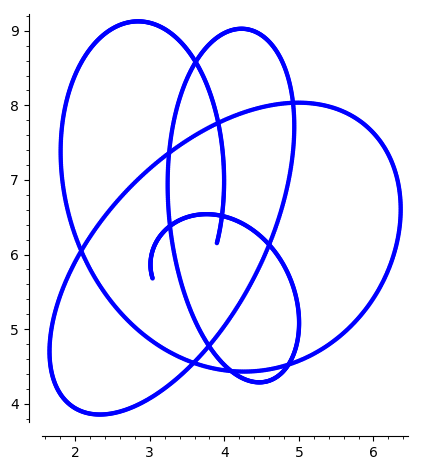}
\caption{$\delta=0.02$}
\end{subfigure}
\caption{Numerically solved trajectories where $c_0=-1$ and the initial data is $(b_{11},b_{12},s,\theta_0,\theta)=(0.65,0,\mu_0+\delta,0,0)$, where $\delta$ varies between $0\dots 0.02$. The solution blows up quickly when $\delta<0$.}
\label{24radatb}
\end{figure}

Unfortunately this equilibrium point is not hyperbolic so the linearization is inconclusive in regard to the stability of the equilibrium. However, the Jacobian has one zero eigenvalue and four purely imaginary eigenvalues if $c_0<0$ and $\alpha<c_1<1$, where $\alpha\approx 0.71$ is the smaller positive root of $q=y^8 - 12y^4 + 3$. So one might suspect that the flow is ''stable'' for $\alpha<c_1<1$ and ''unstable'' for $c_1<\alpha$. Indeed the numerical computations seem to suggest this, see  Figures \ref{24radat} and \ref{24radatb}. We have varied $s$ in these examples; varying other initial values gives similar results.

It is of interest to check whether the Gerstner type solution \eqref{24ell} satisfies the free boundary condition for some choice of $f^1$ and $f^2$. Assuming that $\beta(s)$, $s\in I\subset \mathbb{R}$, is a regular curve in the $z$ plane such that $\partial D=\beta(I)$, pressure should be constant along $\beta$; in other words, we need to satisfy
\begin{equation}
    \langle \nabla p(\beta(s),t),\beta'\rangle =0
    \label{free}
\end{equation}
for all $t$ and $s$. Unfortunately it turns out that this condition cannot be met as we will show. The type of Boussinesq approximation we use does not matter so let us calculate the partial derivatives of $p$ from the standard Boussinesq approximation \eqref{newton2}. Putting $\gamma=c_1^2-1/c_1^2$ we obtain
\begin{align*}
    p_{10}/(\overline p\mu_0^2)=&
\gamma(f^1f_{10}^1 - f^2f_{01}^1)\cos^2(\mu_0t)
-\gamma(f^1f_{01}^1+f^2f_{10}^1)\cos(\mu_0t)\sin(\mu_0t)\\
 & +(c_1^2f^1 - \gamma z_2f_{01}^1  )\cos(\mu_0t)
 - (  \gamma z_2f_{10}^1 + c_1^2f^2)\sin(\mu_0t)
 + c_1^2f^2f_{01}^1 + f^1f_{10}^1/c_1^2\ ,\\
  p_{01}/(\overline p\mu_0^2)=&  \gamma(f^1f_{01}^1 + f^2f_{10}^1)\cos^2(\mu_0t) 
+\gamma(f^1f_{10}^1 - f^2f_{01}^1)\cos(\mu_0t)\sin(\mu_0t)\\
&+(\gamma z_2f_{10}^1 + f^2/c_1^2)\cos(\mu_0t) 
+ (f^1/c_1^2-\gamma z_2f_{01}^1 )\sin(\mu_0t) 
+f^1f_{01}^1/c_1^2 - c_1^2f^2f_{10}^1 - \gamma z_2\ .
\end{align*}
Thus $\nabla p/(\overline p \gamma\mu_0^2)$ contains the expression
\[
G(z)a(t)=\begin{pmatrix}
    g_1 & -g_2 \\
    g_2 & g_1
\end{pmatrix} \begin{pmatrix}
    \cos^2(\mu_0t) \\
    \cos(\mu_0t)\sin(\mu_0t)
\end{pmatrix}\ ,
\]
where $g_1=f^1f_{10}^1 - f^2f_{01}^1$ and $g_2=f^1f_{01}^1+f^2f_{10}^1$. Now the free surface condition \eqref{free} implies $G^T\beta' =0$ for all $s\in I$. Thus 
\[
\det(G)=g_1^2+g_2^2=|\nabla f^1|^2\big((f^1)^2+(f^2)^2\big)=0\ .
\]
This implies that $f^1$ and $f^2$ are constant in $\beta(I)$ and thus constant everywhere since they are an anti-CR pair. This leads to a trivial solution as we further obtain $f^1=f^2=0$.

\subsubsection{Case 2}
Suppose that we have the constraints $g_{13}=g_{24}=0$, or that $v=\big(z_1,z_2,f^1(z_1),f^2(z_2)\big)$. Reducing as in \cite[Lemma 5.4]{toju}, we may suppose at the outset that
\[
A=\begin{pmatrix} 
\cos(\theta) & -\sin(\theta) \\
\sin(\theta) & \cos(\theta)
\end{pmatrix} \begin{pmatrix}
b_{11} &b_{12} & \ell b_{12} & b_{11}/\ell \\
0 & 1/b_{11} & \ell/b_{11} & 0
\end{pmatrix}\ .
\]
This  satisfies the incompressibility condition \eqref{det2dja3d} with
\[
  \det(d\varphi)=1-(f^1)'(f^2)'\neq 0 \ .
\]
Again for density we choose the simplest situation $\rho=c_0z_2$ and moreover we choose $\theta=0$. Inspecting condition \eqref{h2dja3d}, we first find that $b_{12}$ and $b_{11}/\ell$ are linearly dependent, and by a linear transformation we may assume that $b_{12}=0$. Then we obtain for $\ell$ and $b_{11}$
\begin{align*}
    (\ell'/\ell^2)b_{11}^2&=c_{14}\\
    -\ell'/b_{11}^2+c_0\int \ell/b_{11}\textrm{ }dt&=c_{23}\ .
\end{align*}
Further let $b_{11}=y'$. This can be transformed to the system
\begin{align*}
    \ell'=&(c_0/2)y\ell \\
    y'=&\frac{ k_0 \sqrt{y}+c_0y^{3}}{10 y}\ .
\end{align*}

If $k_0=0$ there is an explicit solution
\[
 \ell = \frac{k_3}{\left(c_0 t+k_2 \right)^{5}}\quad,\quad
y = k_1-
\frac{10}{k_2 +c_0t}\ .
\]
Here $k_1$ does not appear in $A$, $k_2=0$ can be assumed by translation of $t$ and by scaling we may take $k_3=c_0^5$. Thus this solution can be written as
\[
A=
\begin{pmatrix}
    10/c_0&0\\0&c_0/10
\end{pmatrix}
\begin{pmatrix}
t^{-2} & 0 & 0 &  t^3 \\
0 & t^2 &  t^{-3} & 0
\end{pmatrix}\ .
\]
This solution is unstably stratified regardless of how $c_0$ is chosen.

\subsubsection{Case 3}
Now we consider the case $g_{13}+g_{24}=g_{14}=0$, where $v=\big(z_1,z_2,z_2(f^1)'(z_1)+f^2(z_1),f^1(z_1)\big)$, and we see from \cite[Lemma 5.5]{toju} that
\[
A=\begin{pmatrix} 
\cos(\theta) & -\sin(\theta) \\
\sin(\theta) & \cos(\theta)
\end{pmatrix}\begin{pmatrix}
    \ell b_{12} & b_{12} & b_{13} & \ell b_{13} \\
    -\ell/b_{13} & -1/b_{13} & 0 & 0
\end{pmatrix}
\]
is the general form that satisfies condition \eqref{det2dja3d} with
\[
  \det(d\varphi)=-z_2(f^1)''-(f^2)'\neq 0 \ .
\]
Let us suppose that $\theta=0$ and that $\rho = c_1z_1+c_2z_2$. Then we get the following system from \eqref{h2dja3d}:
\begin{align*}
    \ell' = &\frac{k_{0}}{b_{13}^{2}}\\ b_{13}' = &
\frac{\left(c_{2} \ell  -c_{1}\right) b_{13}^{4}}{4 k_{0}}\\
b_{12}' =& 
\frac{ \left(c_{2} \ell -c_{1}\right) b_{12}b_{13}^{3}}{4 k_{0}}\ .
\end{align*}
Eliminating $\ell$ from the second equation we obtain
\[
4b_{13}b_{13}''-16(b_{13}')^2+c_2b_{13}^3=0\ .
\]
If $c_2=0$ the above equation can be solved explicitly. Putting moreover $b_{12}=0$ gives the stably stratified solution
\[
A=\begin{pmatrix}
    0 & 0 & t^{-1/3} & (9/20)c_1t^{4/3} \\
    -(9/20)c_1t^2 & -t^{1/3} & 0 & 0
\end{pmatrix}\ .
\]

\subsubsection{Case 4}
Studying the homogeneous situation in \cite{toju}, we dismissed the fourth case as it only provided solutions that reduced to case $m=3$. Here this does not happen; in fact solution \eqref{2d-col} is already an example of this case, albeit with spatial variables named differently.

In this case we can rewrite the spatial constraints as $g_{23}=g_{24}=0$ by changing the order of the components of $v$. In this case we have $v=\big(z_1,z_2,f^1(z_2),f^2(z_2)\big)$, and due to the incompressibility condition \eqref{det2dja3d} $A$ can be assumed to be of the form
\[
A=\begin{pmatrix} 
\cos(\theta) & -\sin(\theta) \\
\sin(\theta) & \cos(\theta)
\end{pmatrix}\begin{pmatrix}
b_1 & b_2 & b_3 & b_4 \\
0 & 1/b_1 & 0 & 0
\end{pmatrix} \ ,
\]
which gives simply $\det(d\fii)=1$. If we choose $\rho=c_0f^2$, \eqref{h2dja3d} yields the conditions
\begin{align*}
    b_1'b_2-b_2'b_1+2\theta'&=c_{12} \\
        b_1'b_3-b_3'b_1&=c_{13}\\
        b_1'b_4-b_4'b_1-c_0y_1 &= c_{14}\ ,
\end{align*}
where $y_1'=a_{21}=b_1\sin(\theta)$. We can choose for example $b_1$ and $\theta$ arbitrarily and easily solve these equations for $b_2$, $b_3$ and $b_4$, yielding
\begin{align*}
    b_2&=b_1\int\frac{2\theta'-c_{12}}{b_1^2}\textrm{ }dt \\
    b_3&=-b_1\int\frac{c_{13}}{b_1^2}\textrm{ }dt \\
    b_4&=-b_1\int\frac{c_0y_1+c_{14}}{b_1^2}\textrm{ }dt\ .
\end{align*}
We may assume $c_{12}=0$ without loss of generality.

\section{3-dimensional case }
\subsection{$m=3$}
Now we turn to the 3D case and once again we take the simplest case to be our first example and consider the case $m=3$. Thus we have $v=(z_1,z_2,z_3)$. We apply the QR-decomposition to $A$ and write $A=RB$, where $R\in\mathbb{SO}(3)$ is a rotation matrix and $B$ is an upper triangle matrix. By the incompressibility condition \eqref{det2dja3d} we then have
\[
\det(d\fii)=\det(B)=b_{11}b_{22}b_{33}=1.
\]
For condition \eqref{h2dja3d} we have, after substituting $b_{33}=1/(b_{11}b_{22})$,
\begin{equation}
\begin{aligned}
    h^1 &= b_{12}'b_{13}-b_{13}'b_{12}+b_{22}'b_{23}-b_{23}'b_{22}+\frac{w_1}{b_{11}}-\frac{w_2b_{12}}{b_{11}b_{22}}+w_3(b_{12}b_{23}-b_{13}b_{22})+y_3\rho_{010}-y_2\rho_{001} \\
    h^2 &= -b_{11}'b_{13}+b_{13}'b_{11}+w_2/b_{22}-w_3b_{11}b_{23}+y_1\rho_{001}-y_3\rho_{100} \\
    h^3 &= b_{11}'b_{12}-b_{12}'b_{11}+w_3b_{11}b_{22}+y_2\rho_{100}-y_1\rho_{010}\ ,
\end{aligned}
\label{3x3}
\end{equation}
where
\[
w=2\big(\langle R_2',R_3\rangle ,-\langle R_1',R_3\rangle ,\langle R_1',R_2\rangle \big)\ .
\]
As in the first 2D case \eqref{2x2}, we may use linear transformations to assume without loss of generality that either $\rho=f(z_3)+c_0z_2$, in which case $a_{31}=a_{32}=0$, or $\rho=c_0z_3$.

Suppose first that $\rho=f(z_3)+c_0z_2$. Since $a_{31}=a_{32}=0$, the rotation matrix $R$ must be of the form
\begin{equation}
R=\begin{pmatrix}
    \cos(\theta) & -\sin(\theta) & 0 \\
    \sin(\theta) & \cos(\theta) & 0 \\
    0 & 0 & 1
\end{pmatrix}
\label{2drot}
\end{equation}
for some function $\theta(t)$, and $y_3'=b_{33}$. In this case $w=(0,0,2\theta')$ and $h$ can be written as
\begin{align*}
    h^1 &= b_{12}'b_{13}-b_{13}'b_{12}+b_{22}'b_{23}-b_{23}'b_{22}+2\theta'(b_{12}b_{23}-b_{13}b_{22})+c_0y_3=c_{23} \\
    h^2 &= -b_{11}'b_{13}+b_{13}'b_{11}-2\theta'b_{11}b_{23}=c_{13} \\
    h^3 &= b_{11}'b_{12}-b_{12}'b_{11}+2\theta'b_{11}b_{22}=c_{12}\ .
\end{align*}
There are different ways to find partial solutions to this system, but one fairly general solution can be found by assuming only that $c_{12}=0$. We may then take $b_{11}$, $b_{22}$, and $\theta$ as arbitrary functions and solve for $b_{12}$, $b_{23}$, and $b_{13}$:
\begin{align*}
    b_{12} &= b_{11} \int \frac{2\theta'b_{22}}{b_{11}}\textrm{ }dt \\
    b_{23} &= b_{22}  \int \frac{-c_{13}b_{12}/b_{11}+c_0y_3-c_{23}}{b_{22}^2}\textrm{ }dt \\
    b_{13} &= b_{11}\int \frac{2\theta'b_{11}b_{23}+c_{13}}{b_{11}^2} \textrm{ }dt\ .
\end{align*}

Then suppose that $\rho=c_0z_3$ is linear. Once again we present the solution to \eqref{3x3} in case $h^3=0$. We choose $b_{11}$ and $b_{22}$, as well as the whole rotation matrix $R$ arbitrarily, and solve for $b_{12}$, $b_{23}$, and $b_{13}$ to obtain

\begin{align*}
    b_{12} &= b_{11} \int \frac{w_3b_{22}}{b_{11}}\textrm{ }dt \\
    b_{23} &= b_{22}  \int \frac{(w_1+b_{12}(c_0y_1-c_{13}))/b_{11}-c_0y_2-c_{23}}{b_{22}^2}\textrm{ }dt \\
    b_{13} &= b_{11}\int \frac{-w_2/b_{22} + w_3b_{11}b_{23}-c_0y_1+c_{13}}{b_{11}^2} \textrm{ }dt\ .
\end{align*}

Finally, if $\rho=c_0z_3$ and we choose $b_{13}=b_{23}=0$ and $R$ as in \eqref{2drot}, we have a flow of the form \eqref{columnar}. Thus we can add two terms to $\fii_3$ as shown in Section 3 and obtain $\fii=Av$, where
\begin{align*}
A&=\begin{pmatrix}
    \cos(\theta) & -\sin(\theta) & 0 \\
    \sin(\theta) & \cos(\theta) & 0 \\
    0 & 0 & 1
\end{pmatrix}\begin{pmatrix}
    b_{11} & b_{12} & 0 & 0 & 0 \\
    0 & b_{22} & 0 & 0 & 0 \\
    0 & 0 & 1/(b_{11}b_{22}) & a_1 & a_2
\end{pmatrix}\ ,\\
v&=\big(z_1,z_2,z_3,f^1(z_1,z_2),f^2(z_1,z_2)\big)\ ,
\end{align*}
where again $b_{11}$, $b_{22}$, and $\theta$ are arbitrary, 
\[
    b_{12} = b_{11} \int \frac{w_3b_{22}}{b_{11}}\textrm{ }dt 
\]
and $a_i$ are linearly independent solutions of \eqref{q-yht} with $a=1/(b_{11}b_{22})$.

\subsection{$m=5$}

Much like for the 2D cases with $m=4$, we will only find a rather small subset of the solutions in the 3D case when $m\geq 5$.
The three cases considered in \cite[Section 5]{toju2} readily yield solutions to the Euler--Boussinesq equations as well. For example in \cite[Section 5.1]{toju2}  there was the solution $\fii=Av$ with
\begin{align*}
A&=\begin{pmatrix}
\cos(\theta)/\sqrt{\theta'} & -\sin(\theta)/\sqrt{\theta'} & 0 & \cos(\theta)/\sqrt{\theta'} & \sin(\theta)/\sqrt{\theta'} \\
\sin(\theta)/\sqrt{\theta'} & \cos(\theta)/\sqrt{\theta'} & 0 & -\sin(\theta)/\sqrt{\theta'} & \cos(\theta)/\sqrt{\theta'} \\
0 & 0 & \theta' & 0 & 0
\end{pmatrix}\ ,\\
v&=\big(z_1,z_2,z_3,f^1(z_1,z_2,z_3),f^2(z_1,z_2,z_3)\big)\ , \\
\det(d\fii)&=1-  \big(f^1_{100}\big)^2-\big(f^1_{010}\big)^2\ne 0\ ,
\end{align*}
where $\theta(t)$ is arbitrary and $f^1$ and $f^2$ are an anti-CR pair with respect to $z_1$ and $z_2$. Recall that a solution to the Euler equations is also a solution to the Euler--Boussinesq equations if $\nabla \rho \times \nabla \fii_3 = 0$. Here $\fii_3=\theta'z_3$, so choosing $\rho$ to be a function of $z_3$ we obtain a solution to the Euler--Boussinesq equations. The other two cases in \cite[Section 5]{toju2}  are
\begin{align*}
  A&=\begin{pmatrix}
  b&0&0&0&b\ell \\
  0&b\ell &0&b&0\\
  0&0&1/(b^2\ell )&0&0
  \end{pmatrix}\ ,\\
v&=\big(z_1,z_2,z_3,f^1(z_1,z_3),f^2(z_2,z_3)\big)\ ,\\
\det(d\fii)&=1-f^1_{100}f^2_{010}\ne 0\ ,
\end{align*}
where $b^2\ell '=1$; and
\begin{align*}
A&=\begin{pmatrix}
0&0&0&b&b\ell \\
b\ell &b&0&0&0 \\
0&0&1/b^2&0&0
\end{pmatrix}\ , \\
v&=\big(z_1,z_2,z_3,f^1(z_1,z_3)+z_2f^2_{100}(z_1,z_3),f^2(z_1,z_3)\big)\ ,\\
\det(d\fii)&=f_{100}^1+z_2f_{200}^2\ne 0\ ,
\end{align*}
where $b^2\ell '=1$. These are also solutions to the Euler--Boussinesq equations with $\rho=\rho(z_3)$.

We can also find the similar type of solutions in the other three cases mentioned at the start of \cite[Section 5]{toju2}  if we require that $\fii_3=a_{33}z_3$. The derivation of the following formulas is very similar to the above three cases, which were derived in \cite{toju2}, so we present them without proofs. In the following formulas $\theta$ is an arbitrary function with $\theta'>0$. If $G_{14}+G_{25}=G_{15}-G_{24}=0$ and $\fii_3=a_{33}z_3$, then the solutions of the Euler equations, and the solutions of the Euler--Boussinesq equations with $\rho=\rho(z_3)$, are
\begin{equation}
\begin{aligned}
A&=\begin{pmatrix}
\cos(k_1\theta)/\sqrt{\theta'} & -\sin(k_1\theta)/\sqrt{\theta'} & 0 & \cos(k_2\theta)/\sqrt{\theta'} & -\sin(k_2\theta)/\sqrt{\theta'} \\
\sin(k_1\theta)/\sqrt{\theta'} & \cos(k_1\theta)/\sqrt{\theta'} & 0 & \sin(k_2\theta)/\sqrt{\theta'} & \cos(k_2\theta)/\sqrt{\theta'} \\
0 & 0 & \theta' & 0 & 0
\end{pmatrix}\ ,\\
v&=\big(z_1,z_2,z_3,f^1(z_1,z_2),f^2(z_1,z_2)\big)\ ,\\
\det(d\fii)&=1-\lvert \nabla f^1\rvert ^2 \neq 0\ ,
\end{aligned}
\label{ell-m=5}    
\end{equation}
where  $f^1$ and $f^2$ are an anti-CR pair.

If $G_{14}=G_{25}=0$, we have
\begin{equation}
\begin{aligned}
A&=\begin{pmatrix}
\cos(\theta) & -\sin(\theta) & 0 \\
\sin(\theta) & \cos(\theta) & 0 \\
0 & 0 & 1
\end{pmatrix} \begin{pmatrix}
e^{\theta}/\sqrt{\theta'} & 0 & 0 & 0 & e^{-\theta}/\sqrt{\theta'} \\
0 & e^{-\theta}/\sqrt{\theta'} & 0 & e^{\theta} /\sqrt{\theta'} & 0 \\
0 & 0 & \theta' & 0 & 0
\end{pmatrix}\ ,\\
v&=\big(z_1,z_2,z_3,f^1(z_1),f^2(z_2)\big)\ ,\\
\det(d\fii)&=1-(f^1)'(f^2)' \neq 0\ .
\end{aligned}
\label{hyp-m=5}    
\end{equation}

And if $G_{14}+G_{25}=G_{15}=0$, then we have
\begin{equation}
\begin{aligned}
A&=\begin{pmatrix}
\cos(\theta) & -\sin(\theta) & 0 \\
\sin(\theta) & \cos(\theta) & 0 \\
0 & 0 & 1
\end{pmatrix} \begin{pmatrix}
0&0&0&1/\sqrt{\theta'}&\theta/\sqrt{\theta'} \\
\theta/\sqrt{\theta'} &1/\sqrt{\theta'}&0&0&0 \\
0&0&\theta'&0&0
\end{pmatrix}\ , \\
v&=\big(z_1,z_2,z_3,f^1(z_1)+z_2(f^2)'(z_1),f^2(z_1)\big)\ ,\\
\det(d\fii)&=(f^1)'+z_2(f^2)''\neq 0\ .
\end{aligned}
\label{par-m=5}    
\end{equation}

In the above three cases $\fii$ is of the form \eqref{columnar} so we can extend these solutions as shown in Section 3; for example the solution  \eqref{ell-m=5}  can be extended as follows: $\fii=Av$, where
\[
A=\begin{pmatrix}
\cos(k_1\theta)/\sqrt{\theta'} & -\sin(k_1\theta)/\sqrt{\theta'} & 0 & \cos(k_2\theta)/\sqrt{\theta'} & -\sin(k_2\theta)/\sqrt{\theta'} & 0 & 0\\
\sin(k_1\theta)/\sqrt{\theta'} & \cos(k_1\theta)/\sqrt{\theta'} & 0 & \sin(k_2\theta)/\sqrt{\theta'} & \cos(k_2\theta)/\sqrt{\theta'} & 0 & 0 \\
0 & 0 & \theta' & 0 & 0 & a_1 & a_2
\end{pmatrix}
\]
and $v=\big(z_1,z_2,z_3,f^1(z_1,z_2),f^2(z_1,z_2),f^3(z_1,z_2),f^4(z_1,z_2)\big)$, where $f^1$ and $f^2$ are an anti-CR pair, and $f^3$ and $f^4$ are arbitrary. Here   $a_i$ are linearly independent solutions of \eqref{q-yht} with $a=\theta'$.

We can add the same expression to \eqref{hyp-m=5} and \eqref{par-m=5} as well.

\subsection{$m=6$, case 1}
\begin{figure}
\centering
\includegraphics[width=0.6\linewidth]{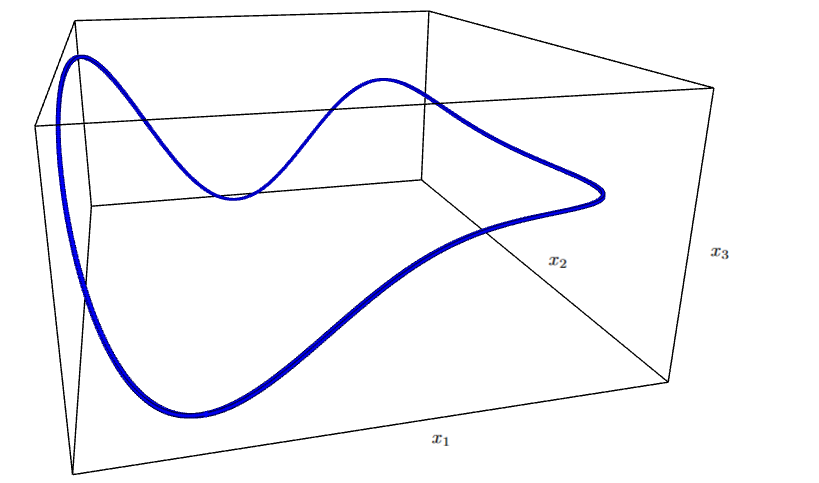}
\caption{A typical particle trajectory for solution \eqref{36ex}.}
\label{36rata}
\end{figure}

Let $v=\big(z_1,z_2,z_3,f^1(z_3),f^2(z_1,z_2),f^3(z_1,z_2)\big)$, where $f^1$ is arbitrary and $f^2$ and $f^3$ are an anti-CR pair. We studied this case for the homogeneous Euler equations in \cite[Section 7]{toju2} and a large part of the calculations we need in the present case was already shown there. We choose $\rho=\sum_ic_iv^i$ and look for solutions of the form
\[
A=\begin{pmatrix}
\cos(\theta) & -\sin(\theta) & 0 \\
\sin(\theta) & \cos(\theta) & 0 \\
0 & 0 & 1
\end{pmatrix} \begin{pmatrix}
    b_{11} & b_{12} &\ell_1b_{11}+\ell_2b_{12} & -\ell_2b_{11}+\ell_1 b_{12} &  \ell_1b_{15} & -\ell_2 b_{15}\\
    0 & b_{22} & \ell_2 b_{22} & \ell_1 b_{22} & \ell_1b_{25}&  -\ell_2b_{25} \\
    0 & 0 &0 & 0 &\ell_1 b_{35}& -\ell_2 b_{35}
\end{pmatrix}\ ,
\]
where the conservation of volume requires that
\begin{equation}
b_{11}b_{22}b_{35}\big(\ell_1^2+\ell_2^2\big)+1=0\ .
\label{36det}
\end{equation}
In this case
\[
\det(d\fii)=f^2_{100}+f^2_{010}f^1_{001}\neq 0\ .
\]
We collect the rest of the constraints from the expressions of $h^j$, which are
\begin{align*}
h^1=&Q_{23}+Q_{24}(f^1)'-(Q_{35}+c_3y_5)f^2_{01}-(Q_{45}+c_4y_5)(f^1)'f^2_{01}+(Q_{36}+c_3y_6)f^2_{10}+(Q_{46}+c_4y_6)(f^1)'f^2_{10} \\
    h^2=&-Q_{13}-Q_{14}(f^1)'+(Q_{35}+c_3y_5)f^2_{10}+(Q_{45}+c_4y_5)(f^1)'f^2_{10}+(Q_{36}+c_3y_6)f^2_{01}+(Q_{46}+c_4y_6)(f^1)'f^2_{01} \\
    h^3=&Q_{12}+(Q_{15}-Q_{26}+c_1y_5-c_2y_6)f^2_{01}-(Q_{16}+Q_{25}+c_1y_6+c_2y_5)f^2_{10}-(Q_{56}+c_5y_6-c_6y_5)\lvert\nabla f^2\rvert ^2\ .
\end{align*}
All the equations that were used to prove \cite[Lemma 7.4]{toju2} are also present here. In that Lemma we showed that from these equations it follows that
\begin{align*}
   & \ell_1=k_1+\cos(\mu)\ , \quad \ell_2=\sin(\mu)\ ,\\
    &  \mu'b_{11}b_{22}=\sqrt{c_{12}^2-c_{13}^2-(k_1c_{12}+c_{23})^2}=:1/k_2
\end{align*}
for some constant $k_1$ and function $\mu$. We also have
\begin{align*}
    \mu'b_{11}^2&=(k_1c_{12}+c_{23})\cos(\mu)+c_{13}\sin(\mu)-c_{12}\\
    b_{12}&=k_2\Big((k_1c_{12}+c_{23})\sin(\mu)-c_{13}\cos(\mu)\Big)b_{22}\, .
\end{align*}
We further assume that $k_1=0$. Now condition \eqref{36det} gives $b_{35}=-k_2\mu'$. Then $a_{35}=-k_2\mu'\cos(\mu)$ and $a_{36}=k_2\mu'\sin(\mu)$, which we can integrate to obtain $y_5=-k_2\sin(\mu)$, $y_6=-k_2\cos(\mu)$. Then we find that we must have $c_1=c_2=c_3=c_4=0$, in which case all equations except the one containing $Q_{56}$ are the same as in the case of Euler equations. Then, as is shown in \cite{toju2}, we have $b_{15}=b_{25}=0$. The only remaining equations after this are those that contain $Q_{12}$ and $Q_{56}$. From them we solve for $\mu$ and $\theta$:
\begin{align*}
    k_2^2(\mu')^3 &= k_2c_5\cos(\mu) - k_2c_6\sin(\mu) + c_{56} \\
    2k_2\theta'&=-\mu' / (c_{23}\cos(\mu) + c_{13}\sin(\mu) - c_{12})=-1/b_{11}^2\ .
\end{align*}
Thus we have a semi-explicit periodic solution where one ODE needs to be solved numerically. As $\mu'$ must not be zero, $\lvert c_{56}\rvert $ has to be large enough for the solution to exist for all $t$. The simplest example is obtained by choosing $c_{13}=c_{23}=0$, $c_{12}=-1$, which gives
\begin{equation}
A=\begin{pmatrix}
    \cos(\theta)/\sqrt{\theta'} & -\sin(\theta)/\sqrt{\theta'} & \cos(\theta)/\sqrt{\theta'} & \sin(\theta)/\sqrt{\theta'} & 0 & 0 \\
    \sin(\theta)/\sqrt{\theta'} & \cos(\theta)/\sqrt{\theta'} & -\sin(\theta)/\sqrt{\theta'} & \cos(\theta)/\sqrt{\theta'} & 0 & 0 \\
    0 & 0 & 0 & 0 & 2\theta'\cos(2\theta) & 2\theta'\sin(2\theta)
\end{pmatrix}\ ,
\label{36ex}
\end{equation}
where
\[
8(\theta')^3+c_5\cos(2\theta)+c_6\sin(2\theta)+c_{56}=0\ .
\]
One particular periodic particle path in this case is shown in Figure \ref{36rata}. This solution is neither stably nor unstably stratified in particular; instead the denser particles are periodically either at the top or at the bottom of the fluid.

\subsection{  $m=6$, case 2}
Let $v=\big(z_1,z_2,z_3,f^1(z_1),f^2(z_2),f^3(z_3)\big)$. We look for solutions of the form
\[
A=\begin{pmatrix}
    a_1 & 0 & 0 & 0 & 0 & \ell_1a_1 \\
    0 & a_2 & 0 & \ell_2a_2 & 0 & 0 \\
    0 & 0 & a_3 & 0 & \ell_3a_3 & 0
\end{pmatrix}\ ,
\]
where the volume preservation condition requires $a_1a_2a_3=\ell_1\ell_2\ell_3=1$, giving
\[
  \det(d\fii)=1+(f^1)'(f^2)'(f^3)'\neq 0\ .
\]
Now
\[
h=\begin{pmatrix}
    \ell_3'a_3^2(f^2)'(z_2)+y_3\rho_{010}-y_5(f^2)'(z_2)\rho_{001} \\
\ell_1'a_1^2(f^3)'(z_3)-y_3\rho_{100} \\
\ell_2'a_2^2(f^1)'(z_1)+y_5(f^2)'(z_2)\rho_{100}
\end{pmatrix}\ ,
\]
from which we see that generally we must have $\rho=c_1z_3+c_2f^2(z_2)$. Thus $A$ is a solution if the following equations are satisfied:
\begin{equation}
\begin{aligned}
    a_1a_2a_3&=1 &
    \ell_1\ell_2\ell_3&=1 &
    \ell_1'a_1^2 &= -c_{16} \\
    \ell_2'a_2^2 &= -c_{24} &
   \ell_3''a_3+2\ell_3'a_3'+c_2-c_1\ell_3& = 0 \ .
\end{aligned}
 \label{al-sys}   
\end{equation}

Here one function can be given arbitrarily. We can try to give solutions in terms of $\ell_3$ since that is the only function with second-order dependence. We get
\begin{equation}
    a_3=\frac{k_0+\int{\frac{c_1\ell_3-c_2}{2\sqrt{\lvert \ell_3'\rvert }}\textrm{ }dt}}{\sqrt{\lvert \ell_3'\rvert }}
    \label{a3ratk}
\end{equation}
from the last equation of \eqref{al-sys} and then we can further  use the remaining equations of \eqref{al-sys} to solve for the other functions in terms of $a_3$ and $\ell_3$. Eliminating $a_1$, $a_2$, and $\ell_2$, we obtain the equation
\[
\ell_3(\ell_1'/\ell_1)^2+\ell_3'(\ell_1'/\ell_1)+c_{16}c_{24}\ell_3^2a_3^2=0\ ,
\]
which is a quadratic polynomial equation for $\ell_1'/\ell_1$. Assuming that its discriminant $\delta=\ell_3'^2-4c_{16}c_{24}\ell_3^3a_3^2\geq 0$, we can solve $\ell_1'/\ell_1$ and integrate to obtain $\ell_1$ and eventually the rest of the functions. Let $s^2=\delta/(4\ell_3^2)$; then we compute
\begin{equation}
\begin{aligned}
    \ell_1 &=\frac{k_1}{\sqrt{\lvert \ell_3\rvert }} \exp\Big(\int{s(t)\textrm{ }dt}\Big)
    &\ell_2 = 1/(\ell_1 \ell_3) \\
    a_1^2 &= -c_{16}/\ell_1'
    &a_2^2 = -c_{24}/\ell_2'\ .
\end{aligned}
\label{muut-ratk}
\end{equation}
$k_1$ may be scaled to $1$. Formulas \eqref{a3ratk} and \eqref{muut-ratk} give a local solution but unfortunately formula \eqref{a3ratk} is meaningful only when $\ell_3'\neq 0$. An example where $\ell_3$ is not monotone can be found by supposing that $a_3=1$, $c_1=-N^2$ is negative and also $c_{16}c_{24}<0$. Then the final equation of \eqref{al-sys} gives
\begin{equation}
    \ell_3=-c_2/N^2+k_2\cos(Nt)+k_3\sin(Nt)\ .
    \label{l3ratk}
\end{equation}
Choosing the constants such that $\ell_3$ is always positive, \eqref{l3ratk} along with \eqref{muut-ratk} gives a global solution describing a stably stratified fluid.

\end{document}